\documentclass[mathpazo]{cicp}
\usepackage[margin=1in]{geometry}
\usepackage{hyperref}
\usepackage{amsmath,epsf,cite}
\usepackage{amssymb,amsthm,tocvsec2}
\usepackage{graphicx}
\usepackage{epsfig,epsf,latexsym,subfigure}
\usepackage{float}
\usepackage{color}
\usepackage{mathrsfs}
\usepackage{indentfirst}

\usepackage[width=.85\textwidth]{caption}

\usepackage{datetime}
\usepackage{float}
\usepackage{mathtools}

\usepackage{mathdots}
\usepackage{cleveref}
\usepackage{cancel}
\usepackage{lineno}

\numberwithin{equation}{section}

\theoremstyle{plain}

\theoremstyle{definition}                  

\def\A1{\mathcal{A}_1}

\begin{document}
\title{Numerical Approximations of the Allen-Cahn-Ohta-Kawasaki (ACOK) Equation with Modified Physics Informed Neural Networks (PINNs)}
\author[J. Xu, J. Zhao and X. Zhao]{Jingjing Xu \affil{1} and 
Jia Zhao \affil{2} and Yanxiang Zhao \affil{1}\comma \corrauth}
\address{
\affilnum{1}\ Department of Mathematics, George Washington University, Washington, DC, USA;\\
\affilnum{2}\ Department of Mathematics \& Statistics, Utah State University, Logan, UT, USA. }
\email{ {\tt yxzhao@email.gwu.edu} (Y.~Zhao)}


\begin{abstract}
The physics informed neural networks (PINNs) has been widely utilized to numerically approximate PDE problems. While PINNs has achieved good results in producing solutions for many partial differential equations, studies have shown that it does not perform well on phase field models. In this paper, we partially address this issue by introducing a modified physics informed neural networks. In particular, they are used to numerically approximate Allen-Cahn-Ohta-Kawasaki (ACOK) equation with a volume constraint.
Technically, the inverse of Laplacian in the ACOK model presents many challenges to the baseline PINNs. To take the zero-mean condition of the inverse of Laplacian, as well as the volume constraint, into consideration, we also introduce a separate neural network, which takes the second set of sampling points in the approximation process. Numerical results are shown to demonstrate the effectiveness of the modified PINNs. An additional benefit of this research is that the modified PINNs can also be applied to learn more general nonlocal phase-field models, even with an unknown nonlocal kernel.
\end{abstract}

\maketitle

\section{Introduction}

Ohta-Kawasaki(OK) free energy\cite{OK}, usually associated with a volume constraint, has been used to simulate the phase separation of diblock copolymers, i.e., polymers consisting of two types of monomers, $A$ and $B$, that are chemically incompatible and connected by a covalent chemical bond. It is formulated as
\begin{align}
E[u] = & \int_{\mathbb{T}^d} \left[ \dfrac{\epsilon}{2}|\nabla u|^2 + \dfrac{1}{\epsilon}W(u) \right] dx + \dfrac{\gamma}{2} \int_{\mathbb{T}^d} \Big|(-\Delta)^{-\frac{1}{2}}\Big(f(u)-\omega\Big)\Big|^2\ dx \notag\\
&+ \frac{M}{2} \left[\int_{\mathbb{T}^d} (f(u)-\omega) \ dx\right]^2. \label{eq:ACOK}
\end{align}
Here $\mathbb{T}^d = \prod_i [-X_i,X_i) \subset \mathbb{R}^d, d = 1, 2, 3$ is a periodic domain, and $u = u(x)$ is a phase field labeling function representing the density of $A$ species with interfacial width $\epsilon$. Function $W(u) = 18(u^2-u)^2$ is a double well potential. The nonlinear function $f(u) = 6 u^5 - 15 u^4 + 10 u^3$ keeps the $0\text{-}1$ structure of $u$ as $f(0) = 0$ and $f(1) = 1$. More importantly it also has the property that $f'(0) = f'(1) = 0$, which localize the force and avoid possible non-zeros and non-ones near the boundary of the interface \cite{ACOK} when considering $L^2$ gradient flow dynamics of (\ref{eq:ACOK}). The first integral is a local surface energy, the second term represents the long-range interaction between the molecules with $\gamma$ controlling its strength and $\omega \in (0,1)$ being the fraction of species $A$, and the third term is a penalty term to fulfill the volume constraint
\begin{align}
\int_{\Omega} (f(u)-\omega) \ dx = 0.
\end{align}

Allen-Cahn-Ohta-Kawasaki(ACOK) model \cite{ACOK}, which is the focus of this paper, takes the $L^2$ gradient dynamics of the OK free energy in \eqref{eq:ACOK} as
\begin{align}
u_t =& - \dfrac{\delta E[u]}{\delta u} \\
= &\varepsilon\Delta u - \frac{1}{\varepsilon}W'(u) - \gamma(-\Delta)^{-1}(f(u)-\omega)f'(u) - M \left[\int_{\mathbb{T}^d} (f(u)-\omega) \ dx\right]f'(u),
\end{align}

Various successful and efficient PDE solvers can be applied to the model, such as Finite Difference, Finite Element 
\cite{FE}, Spectral method \cite{Orszag, nla.cat-vn4108574}, etc. In recent years, some novel numerical methods have also been studied for Cahn-Hilliard dynamics, i.e., the $H^{-1}$ gradient flow dynamics of OK energy, such as the midpoint spectral method \cite{mid}, and the Invariant Energy Quadratization (IEQ) method \cite{IEQ}. In \cite{ACOK}, a first-order energy stable linear semi-implicit method is introduced and analyzed for the ACOK equation. Our goal in this paper is to study the approximation of the solution of the ACOK equation with neural networks. 

\begin{figure}[H]
\centering
\includegraphics[scale=0.15]{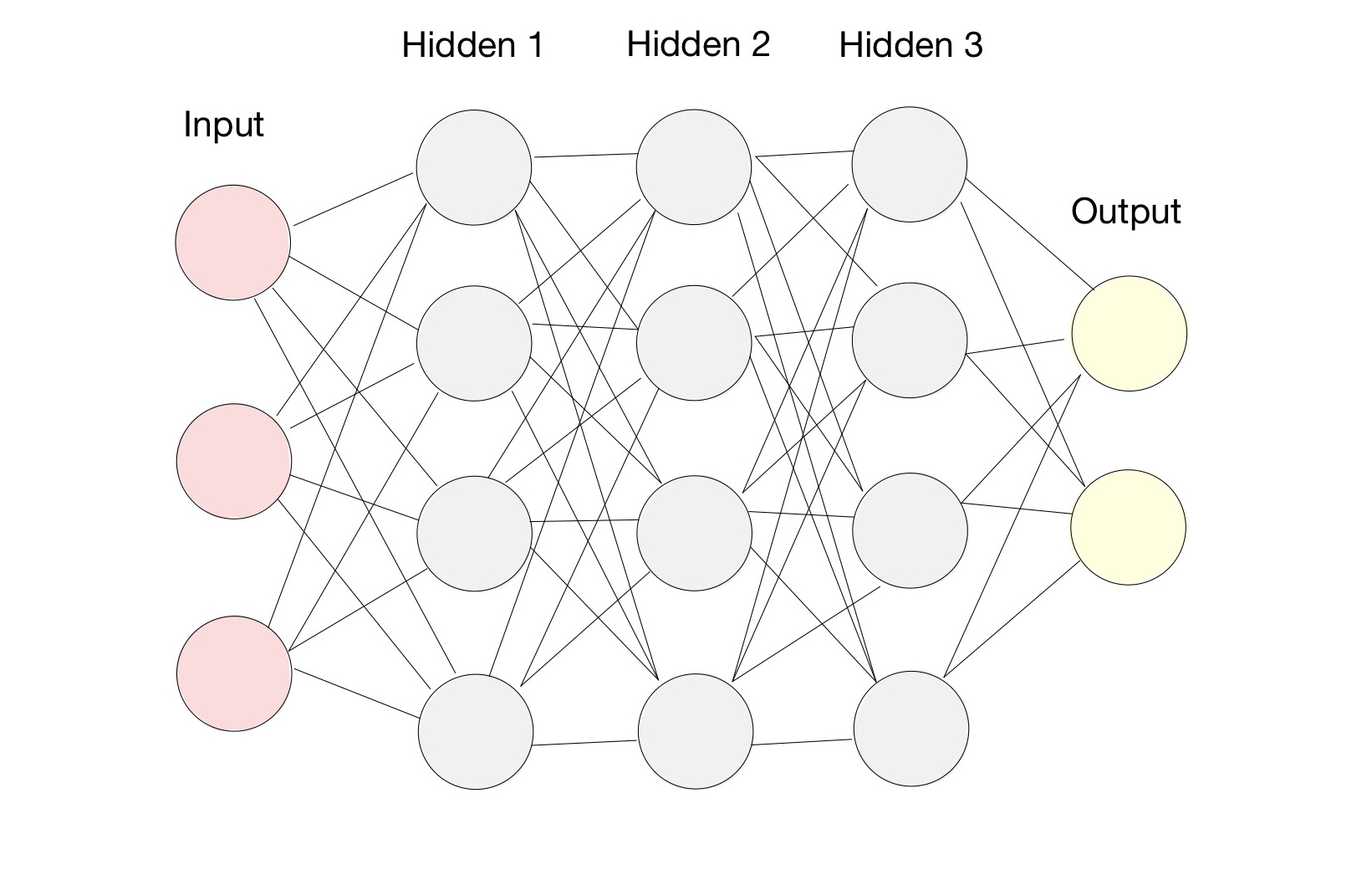}
\caption{Fully-connected DNN structure.}\label{Fig:FC}
\end{figure}

Recently, there has been an increasing interest in solving PDEs with machine learning methods, among which the physics-informed neural networks(PINNs) \cite{PINNs} have gained tremendous popularity.
PINNs are based on a fully-connected feed-forward deep neural network (DNN)\cite{doi:10.1073/pnas.79.8.2554, 10.5555/104279}, which is often interpreted utilizing universal approximation theorem \cite{Cybenko1988b, Kurt1991251, haykin99a, Hassoun:1995}. The structure of a fully-connected DNN, consisting of an input layer, hidden layers, and an output layer, is shown in \Cref{Fig:FC}. Every node in one layer is connected with every node in another by a series of computations. The nodes in the input layer are the specific data studied, and the nodes in the output are the expected outcome. For example, in the case of facial recognition, the input could be the RBG values (features) of each pixel of the sample pictures, while the output could be the assigned numerical values representing the individuals (classes). The nodes in the hidden layers are called neurons, and they are responsible for the performance of the neural network. There can be many hidden layers in a neural network; in the case of Figure \ref{Fig:FC}, the number of hidden layers is three. 

A column vector is fed to each node of the input layer. Then in the first hidden layer, we multiply it by another row vector, i.e., weight, to get the product. This is the case that there is only one neuron in the relevant hidden layer. If there are more neurons, the row vector will be a weight matrix instead, where the number of rows of the matrix represents the number of neurons. A bias vector is then added to the product, and an activation function is applied to the new vector, which contributes to the output of the current hidden layer, which is also the input of the next layer. In a nutshell, let $M$ be the number of hidden layers, $p$ be the input of the neural network, then the output $a$ of the neural network will be $a^M$, with
\begin{align}
a^0  &= p, \\
a^{i} &= f^{i}(W^{i}a^{i-1}+b^{i}), \ i = 1,..., M-1,
\end{align}
where $a^{i-1}$ is the input, $W^i$ represents the weights, $b^i$ represents the bias, and $f^i$ is the activation function in the $i\text{th}$ layer. 

To obtain the proper weights and bias, we perform an optimization process, for example, gradient descent, on the weight and bias to update them. To calculate the derivatives involved in the optimization process, the network goes through a process called back-propagation \cite{Rumelhart:1986we}, which is essentially a chain rule applied on each layer starting from the final one. We start with a certain loss function:
$$\text{Loss}(y, a(\Theta)),$$
where $y$ is the object, the network is approximating, and $\Theta$ denotes the weights and bias of the DNN. The loss function, depending on the specific problem, could be a mean square error (MSE), mean absolute error (MAE), log loss function, etc. We then perform back-propagation on the loss function, trace all the way back to get the derivatives of weights and bias for each layer, and then use the optimization tool to update the weights and bias. A gradient descent update is as follows:
$$\Theta^{n+1} = \Theta^n - \eta \nabla_\Theta \text{Loss},$$
where $\eta$ is the learning rate, and $\nabla_\Theta \text{Loss}$ is the direction of update. It works essentially like a directing system, telling you which direction to go and how far to go in that direction based on your current location and destination. The optimization process will be repeated several times until the prescribed accuracy is achieved. We call this number of repetitions the number of epochs.

Furthermore, there are a variety of optimizers that can be selected. The previously mentioned gradient descent is one of them. While it has many advantages, it does not work well in our non-convex problem due to its notoriously poor performance on escaping local minimums \cite{LocalMin}. Therefore, other optimizers, such as momentum \cite{Momentum}, ADAM \cite{ADAM}, Adadelta \cite{Adadelta}, RMSprop \cite{RMSprop}, etc are developed to better address this issue. The momentum method adds an accelerator in the relevant direction. For example, think about rolling a ball down a ramp. As the ball goes down, it gains momentum and travels faster and faster. Other methods, like Adadelta, use different forms of adaptive learning rates. It is worth mentioning that RMSprop works by dividing ``the learning rate for a weight by a running average of the magnitudes of recent gradients for that weight" \cite{RMSprop}. 

In our case, we adopt ADAM as our optimizer for each epoch and L-BFGS-B with a certain stopping condition as the optimizer afterward. Introduced in \cite{ADAM}, ADAM can be viewed as a combination of the momentum method and RMSprop. It uses the mean $m_t$ and uncentered variance $v_t$ of the gradient with regards to time $t$ to adapt the learning rate for each weight and bias of the neural network:
$$
m_t = \beta_1m_{t-1}+(1-\beta_1)g_t, \quad v_t = \beta_2v_{t-1}+(1-\beta_2)g_t^2,
$$
where $g_t$ is the gradient of the loss function at time step $t$, and the hyper-parameters $\beta_1,\beta_2 \in [0,1)$ direct the decay rates of $m_t$ and $v_t$. A bias correction is also applied to $m_t$ and $v_t$, which leads to
$$\hat{m}_t = \frac{m_t}{1-\beta_1^t}, \quad \hat{v}_t = \frac{v_t}{1-\beta_2^t},$$
and to update the weights and bias, we perform
$\Theta_t = \Theta_{t-1} - \eta \frac{\hat{m}_t}{\sqrt{\hat{v}_t}+\varepsilon}$,
where $\eta$ is the learning rate, and $\varepsilon$ is arbitrary.
On the other hand, the L-BFGS-B optimizer \cite{LBFGSB} uses a ``limited memory" Broyden–Fletcher– Goldfarb–Shanno (BFGS) \cite{10.1093/imamat/6.1.76, 10.1093/comjnl/13.3.317, Goldfarb1970AFO, 10.2307/2004840} matrix to represent the Hessian matrix of the objective function, making it great for optimization problems with a large number of variables. Stems from the L-BFGS method \cite{LBFGS1, LBFGS2}, L-BFGS-B places upper and lower bound constraints on variables. It identifies fixed and free variables for each step with a basic gradient method. It then performs the L-BFGS method on the free variables to achieve better accuracy.

For PINNS, sampled collocation points are fed to the DNN, and the function value is returned to the output layer. Weight and bias in each layer can be learned or updated by optimizing the error function to reach the prescribed accuracy, in which the back-propagation technique is needed to compute the first-order derivatives. Different from traditional machine learning DNN, which requires a very large amount of ground truth data to make an accurate prediction, PINNs take advantage of the physics information and make up for the lack of sufficient ground truth data. It divides the error function into two components, one corresponding to the known ground truth data and the other corresponding to the physics information. Consequently, our goal is to minimize the error between the limited amount of ground truth and the prediction, as well as how much the prediction "fits" the physics information. Furthermore, as shown in \cite{inverse}, PINNs can be employed to solve inverse problems, in particular, inverse scatter problems in photonic metamaterials and nano-optics technologies, by transforming them into parameter retrieving problems. However, PINNs also have their downsides. As suggested in \cite{FailurePINNs}, ``existing PINN methodologies can learn good models for relatively trivial problems", and it can fail when predicting solutions for a more complex model, which is likely due to the fact that the setup of PINNs makes it hard to optimize the less-smooth loss landscape. They have, therefore, proposed a ``curriculum regularization", which is finding a good initialization for weights and training the model gradually. Notably, this has made the optimization process easier. They also mentioned a sequence-to-sequence learning method, which is similar to the Time-Adaptive Method proposed in \cite{timeadap}. This method can make learning easier by learning in small intervals.


While many successful and efficient traditional numerical methods exist, our motivation to study the ACOK equation using machine learning methods is that they can learn the unknown operators in a family of generalized equations from a set of data samples, i.e., the inverse problem. This would benefit future studies focusing on nonlocal PDEs with general nonlocal operators. It can also potentially help tackle real-world challenges in fields such as physics, medicine, fluid dynamics, aerospace, etc. Our contribution to this work is to tackle the challenges that the complexity of the ACOK model has imposed on PINNs. As shown in \ref{eq:ACOK}, the long-range interaction term, the volume constraint, and the zero-mean condition make it necessary to apply computational and structural modifications to the basic PINNs. Our solution for the long-range interaction term is to apply the Laplace operator on a second neural network output and ensure its accuracy. We introduce a separate neural network into the model for the volume constraint and the zero-mean condition, which takes a separate set of uniform sampling points as its input.

The rest of this paper is organized as follows. In \Cref{methods}, we present the numerical methods. Specifically, we briefly revisit the baseline PINNs in \Cref{baseline}. Then we elaborate on the modified PINNs, including the application of the periodic boundary condition, the approximation of the long-range interaction term, the volume constraint, and the zero-mean condition, in \Cref{modified}. We will then present and discuss some 1-dimensional results in \Cref{Results}. In the end, we give a detailed discussion and explain our future work in \Cref{conclusion}.

\newpage
\section{Numerical Methods} \label{methods}

Our goal in this paper is to take the ground truth $u(0,x)$ 
data as input, and predict the solution $u(t,x)$ of ACOK equation with neural networks.

\subsection{Baseline PINNs} \label{baseline}
\begin{figure}[H]
\centering
\includegraphics[scale=0.35]{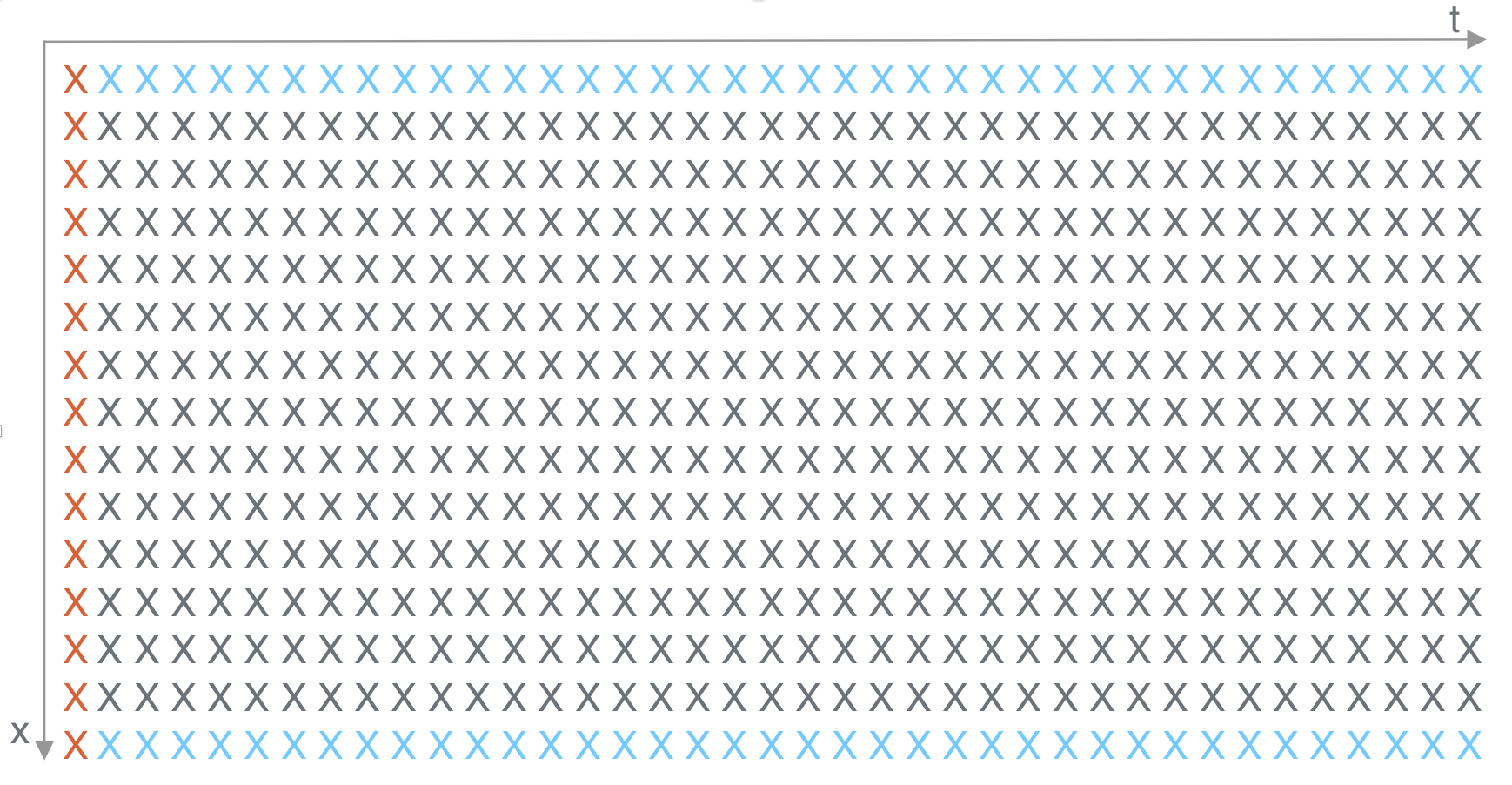}
\caption{Sampling for the initial points, boundary points, and internal points}\label{Fig:IB}
\end{figure}

For the basic setup, we follow the basic PINNs. Let us recall \eqref{eq:ACOK} and define the residual function $F(t,x)$ as
\begin{align}
F(t,x): = &  \ u_t - \varepsilon \Delta u + \frac{1}{\varepsilon} W'(u) + \gamma (-\Delta)^{-1}(f(u)-\omega)\cdot f'(u) \notag \\
& + M \left[\int (f(u)-\omega) \ dx\right] \cdot f'(u),\label{eqn:F}
\end{align}
where $u(t,x)$ can be approximated by a deep neural network $\text{Net}_u$, and $F(t,x)$ can be calculated by a shared parameter neural network $\text{Net}_F$ based on the predicted $u(t,x)$.
The network can be learned by minimizing the mean square error:
$$\text{MSE} = \text{MSE}_u +\text{MSE}_{F},$$
where
$$\text{MSE}_u = \frac{1}{N_u}\sum_{i=1}^{N_u}\left|u(t_u^i,x_u^i)-u^i\right|^2,\quad
\text{MSE}_{F} = \frac{1}{N_F}\sum_{i=1}^{N_F}\left|F(t_F^i,x_F^i)\right|^2.$$
As shown in Figure \ref{Fig:IB}, $\{t_u^i,x_u^i\}_{i=1}^{N_u}$ are selected from the collocation points that are along the boundary, i.e. the blue points, and at the initial time, i.e. the red points. Also, $\{t_F^i,x_F^i\}_{i=1}^{N_F}$ are sampled from the internal collocation points, i.e. the black points. 
For the purpose of computational efficiency and prediction accuracy, a random sample is taken from the initial and boundary collocation points. Moreover, the internal points are selected near-randomly using Latin Hypercube Sampling(LHS) method.

We take pairs of $\{t_u^i,x_u^i\}_{i=1}^{N_u}$ and $\{t_F^i,x_F^i\}_{i=1}^{N_F}$ as the input, and employ the feed forward deep neural network(DNN) to approximate the $u$ value at these collocation points: $u(t_u^i,x_u^i)$ and $u(t_F^i,x_F^i)$. To check the accuracy of the prediction, we take the mean square error between the ground truth $u^i$ and the prediction $u(t_u^i,x_u^i)$. Note that we take the mean square error of $F(t_F^i,x_F^i)$ obtained through $\text{Net}_F$ to ensure that the approximation by the DNN satisfies the physics, which is described by the ACOK equation.

As previously mentioned, ADAM is selected as our optimizer for each epoch, and L-BFGS-B is employed in the last step. Both of them have built-in packages developed in Python for implementation. We also use $\tanh$ as the activation function for each hidden layer. 

\subsection{Modified PINNs}\label{modified}
Next, we demonstrate how we develop the modified PINNs. There are several difficulties in applying the baseline PINNs directly to the ACOK equation, which motivate us to propose the modified PINNs to overcome these issues. 

\subsubsection{Periodic boundary conditions}
As we have the periodic boundary condition in the ACOK equation, $\text{MSE}_u$ is separated further into two components, one corresponding to the initial collocation points, i.e. $\text{MSE}_{in}$, and the other corresponding to the boundary collocation points, i.e. $\text{MSE}_{b}$. Hence, we can update the loss function as:
$$\text{MSE} = \text{MSE}_{in} + \text{MSE}_b +\text{MSE}_{F},$$
where
$$\text{MSE}_{in} = \frac{1}{N_{in}}\sum_{i=1}^{N_{in}}\left|u(t_{in}^i,x_{in}^i)-u_0^i\right|^2,\quad
\text{MSE}_{b} = \frac{1}{N_{b}}\sum_{i=1}^{N_{lb}}\left|u(t_{lb}^i,x_{lb}^i)-u(t_{ub}^i,x_{ub}^i)\right|^2,$$
and
$$\text{MSE}_{F} = \frac{1}{N_F}\sum_{i=1}^{N_F}\left|F(t_F^i,x_F^i)\right|^2.$$
Since we have the ground truth at the initial time, we take the mean square error between the $u$ value predicted by the DNN, $u(t_{in}^i,x_{in}^i)$, and the ground truth $u^i$. However, the periodic boundary condition means that we will calculate $\text{MSE}_{b}$ by taking the mean square error between the $u$ value on the lower bound, i.e. $u(t_{lb}^i,x_{lb}^i)$, and the $u$ value on the upper bound, i.e. $u(t_{ub}^i,x_{ub}^i)$, at the same time $t$, since the ground truth value of $u$ is unknown.

\subsubsection{Approximation of the inverse of Laplacian}
One major difficulty is how to approximate the long-range interaction term in our model. As stated above, the purpose of the shared parameter neural network $\text{Net}_F$ is to approximate $F(t_F^i,x_F^i)$ given $u(t_F^i,x_F^i)$. Recall \eqref{eqn:F}:
$$F: = u_t - \varepsilon \Delta u + \frac{1}{\varepsilon} W'(u) + \gamma (-\Delta)^{-1}(f(u)-\omega)\cdot f'(u) + M \left[\int (f(u)-\omega) \ dx\right] \cdot f'(u).$$
While $u_t$, $\Delta u$ can be obtained by back-propagation, and $W'(u)$, $f(u)$ and $f'(u)$ can be directly calculated using the explicit expression, the difficulty lies with the approximation of $(-\Delta)^{-1}(f(u)-\omega)$ and $\int (f(u)-\omega) \ dx$. 

To approximate the inverse of the Laplacian
$$(-\Delta)^{-1}(f(u)-\omega),$$
we add a second output $\nu$ to $\text{Net}_u$ to approximate
$$\nu := (-\Delta)^{-1}(f(u) - \omega).$$
If a $-\Delta$ operator is applied on both sides, we recover:
$$-\Delta \nu = f(u) - \omega.$$

Since $\nu$ has already been approximated by the DNN, $\Delta \nu$ can be calculated by back-propagation. To ensure the accuracy of $\nu$, a mean square error can be taken between $(-\Delta) \nu$ and $f(u) - \omega$, as the latter is achieved by the explicit expression:
$$\text{MSE}_{Lap} = \frac{1}{N_F} \sum_{i=1}^{N_F}\left|- \Delta \nu (t_F^i,x_F^i) - (f(u(t_F^i,x_F^i)) - \omega)\right|^2,$$
and the first four terms of $F$ can be calculated by:
$$F = u_t - \varepsilon  u_{xx} + \frac{1}{ \varepsilon}  W'(u) + \gamma  \nu \cdot f'(u).$$

Furthermore, as the ground truth of $(-\Delta)^{-1} u_{in}$ and the periodic boundary condition on $(-\Delta)^{-1} u$ are available, the total mean square loss is modified as:
$$\text{MSE} = (\text{MSE}_{u_{in}} + \text{MSE}_{((-\Delta)^{-1}u)_{in}}) + (\text{MSE}_{u_{b}} + \text{MSE}_{((-\Delta)^{-1}u)_{b}}) + \text{MSE}_{F} +\text{MSE}_{Lap},$$
where
$$\text{MSE}_{((-\Delta)^{-1}u)_{in}} = \frac{1}{N_{in}}\sum_{i=1}^{N_{in}}\left|(-\Delta)^{-1}u(t_{in}^i,x_{in}^i)-((-\Delta)^{-1}u)_0^i\right|^2,$$
$$\text{MSE}_{((-\Delta)^{-1}u)_{b}} = \frac{1}{N_{b}}\sum_{i=1}^{N_{lb}}\left|(-\Delta)^{-1}u(t_{lb}^i,x_{lb}^i)-(-\Delta)^{-1}u(t_{ub}^i,x_{ub}^i)\right|^2.$$

\subsubsection{Integral approximation}
\begin{figure}[H]
\begin{center}
\subfigure[LHS Sampling.]{\includegraphics[scale=0.276]{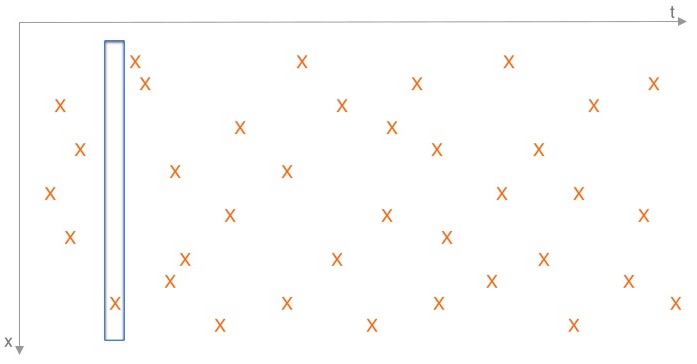}}
\subfigure[Uniform $x$ and Random $t$ Mesh.]{\includegraphics[scale=0.275]{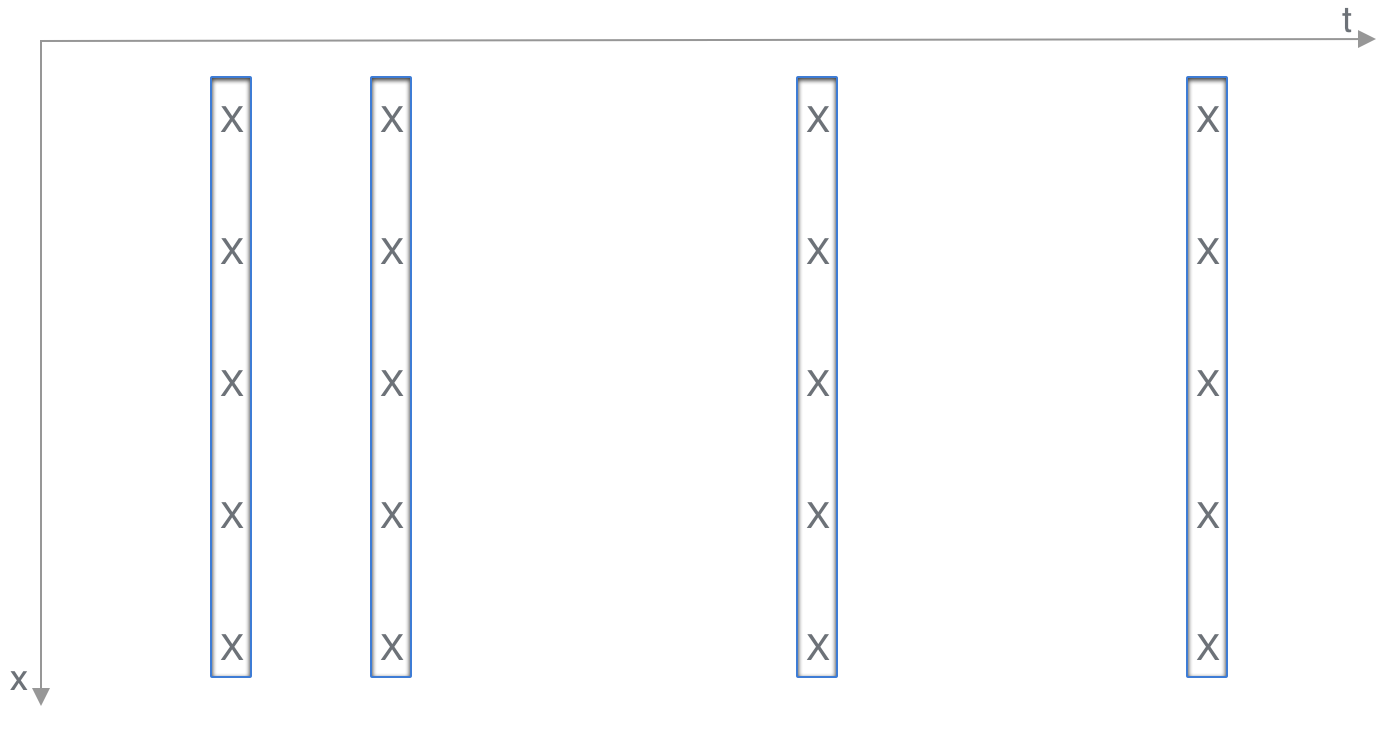}}
\end{center}
\caption{Collocating point sampling strategy: (a) LHS sampling; (b) uniform spatial sampling and random t}
\label{Fig:Uniform2}
\end{figure}

Next, we propose a structural modification to approximate the integral term. The integral term
$$\int (f(u)-\omega) \ dx,$$
imposes a problem to the modified PINNs, as the near-randomly sampled collocation points are not evenly distributed along each $t$. It is usually the case that at some $t$ there are one or no points selected, shown in Figure \ref{Fig:Uniform2}.
As a result, the integral can not be calculated accurately for each $t$. To resolve this issue, a separate shallow neural network $\text{Net}_v$ is introduced, which only takes $t$ as input and has a built-in uniform mesh on $x$. This network aims to output the value of the integral term for each randomly sampled $t$. After the network is trained, the updated optimal weights and bias ensure that for any input $t$, the network will produce the corresponding $\int (f(u)-\omega) \ dx$, even if the new input $t$ is different from the $t$ used for training.
\begin{figure}[H]
\centering
\includegraphics[scale=0.45]{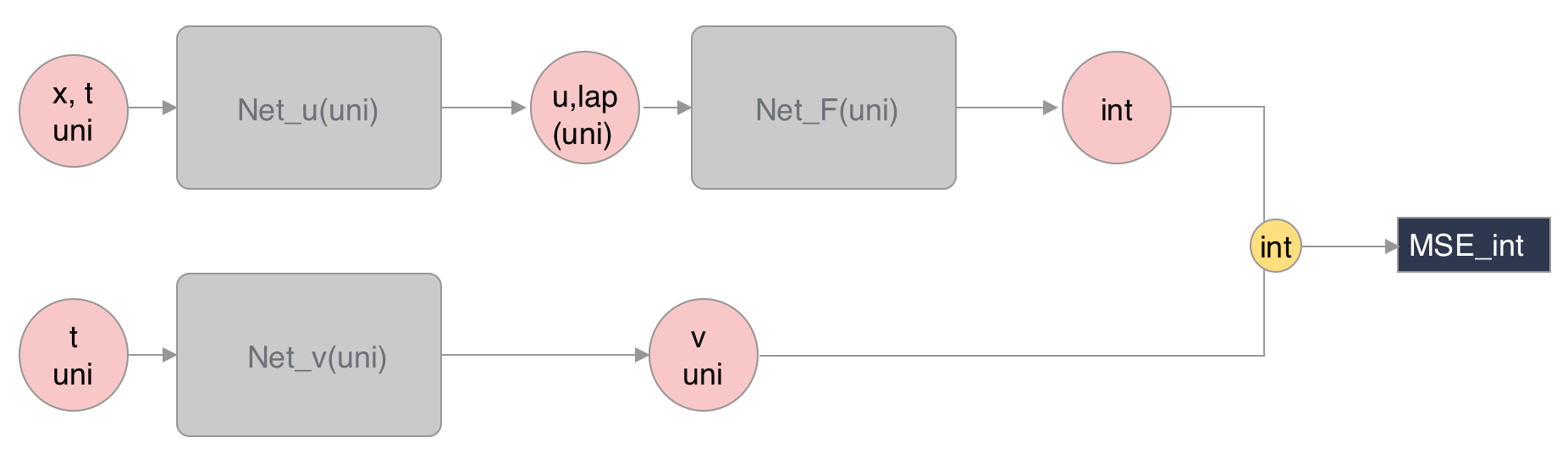}
\caption{$\text{Net}_v$ Structure}\label{Fig:net_v1}
\end{figure}
As demonstrated in Figure \ref{Fig:net_v1}, to train $\text{Net}_v$, the mean square error between the output of $\text{Net}_v$ and the discrete form of the integral
$$\sum_{i=1}^{N_x} f(u(t_{uni}^j,x_{uni}^i)) \cdot \Delta x- \omega,$$
is taken:
$$\text{MSE}_{int} = \frac{1}{N_t} \sum_{j=1}^{N_t}\left|v(t_{uni}^j) - \left[\sum_{i=1}^{N_x} f(u(t_{uni}^j,x_{uni}^i)) \cdot \Delta x- \omega\right] \right|^2,$$
where $\{t_v^j,x_v^i\}$ denotes the new uniform $x$ and random $t$ mesh as in Figure \ref{Fig:Uniform2}, and $u(t_{uni}^j,\linebreak x_{uni}^i)$ can be output by $\text{Net}_u$ with the same set of input. Furthermore, $\sum_{i=1}^{N_x} f(u(t_{uni}^j,x_{uni}^i))$ is obtained by taking the column sum of output of $f(u(t_{uni}^j,x_{uni}^i))$ along each sampled $t$ direction.

After the optimal weight and bias are achieved, the $t$ part of the original LHS sampled points is fed to $\text{Net}_v$, and outputs 
$$v = \int (f(u(t_F^i,x_F^i))-\omega) \ dx,$$
and $F$ can be approximated by
$$F = u_t - \varepsilon  u_{xx} + \frac{1}{\varepsilon}  W'(u) + \gamma  \nu \cdot f'(u) + M v\cdot   f'(u),$$
which is calculated in $\text{Net}_F$.
Taking the mean square error $F$ will obtain the complete $\text{MSE}_F$. 

Therefore, the total mean square loss is modified as:
\begin{align*}
\text{MSE} = & \ (\text{MSE}_{u_{in}} + \text{MSE}_{((-\Delta)^{-1}u)_{in}}) + (\text{MSE}_{u_{b}} + \text{MSE}_{((-\Delta)^{-1}u)_{b}}) \\
& + \text{MSE}_{F}  +\text{MSE}_{Lap} + \text{MSE}_{int},
\end{align*}
where
$$\text{MSE}_{F} = \frac{1}{N_F}\sum_{i=1}^{N_F}\left|F(t_F^i,x_F^i)\right|^2.$$

\subsubsection{Zero-mean constraint}
\begin{figure}[H]
\centering
\includegraphics[scale=0.4]{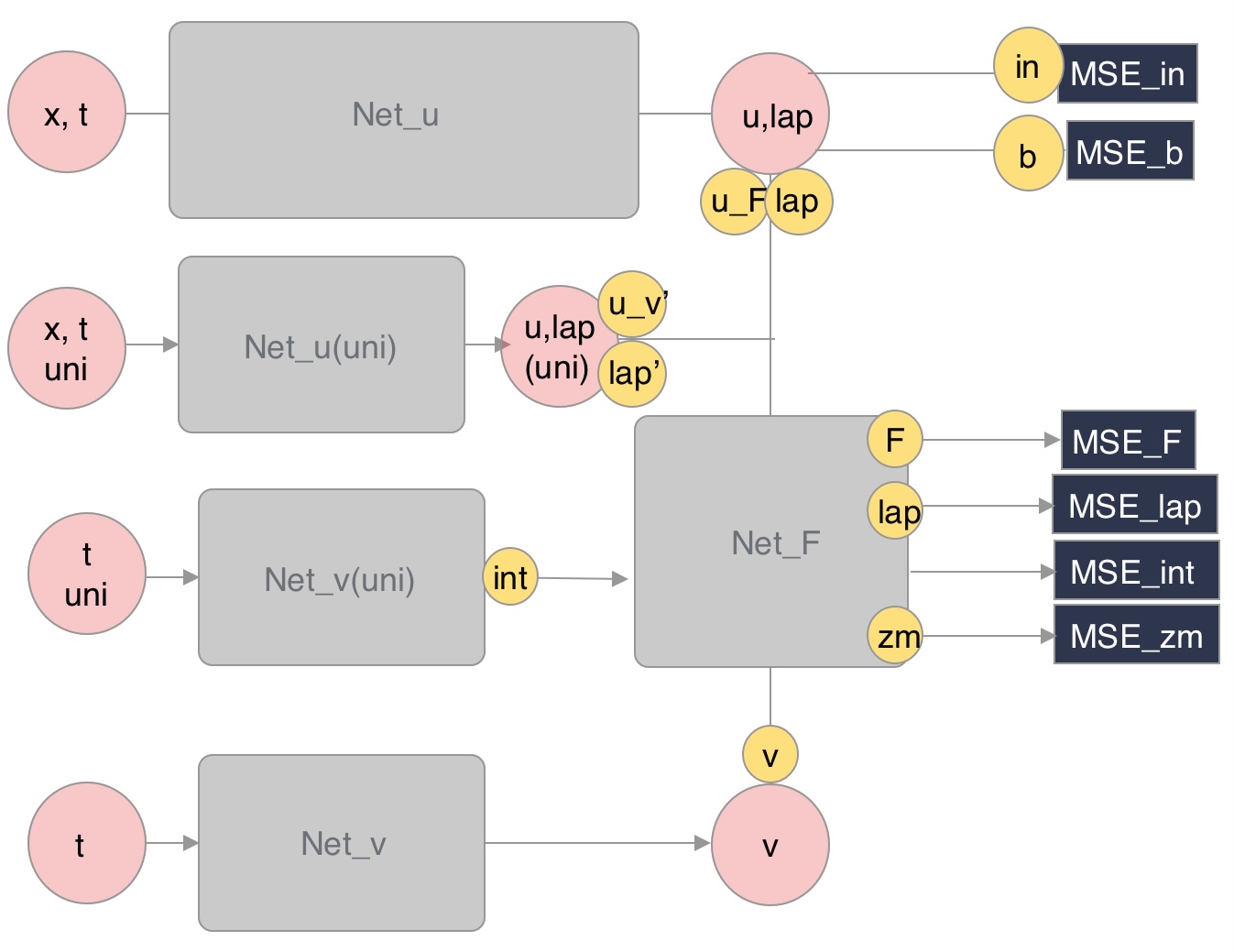}
\caption{The modified PINNs structure with $\text{Net}_v$ and zero-mean condition.}\label{Fig:zm}
\end{figure}

For the zero-mean condition on the inverse of the Laplacian, since the calculation of the mean is involved, the uniform $x$ and random $t$ meshgrid created for the integral approximation can be used as input of $\text{Net}_u$, and the column sum of the output $(-\Delta)^{-1}u(t_F^i,x_F^i)$, which corresponds to each time $t$, is minimized in a loss function 
$$\text{MSE}_{zm} = \frac{1}{N_t} \sum_{j=1}^{N_t}\left|\sum_{i=1}^{N_x} (-\Delta)^{-1}u(t_F^i,x_F^i)\cdot \Delta x\right|^2.$$
The complete modified PINNs structure is shown in \Cref{Fig:zm}.

Hence, the modified PINNs is learned by minimizing the mean square loss:
\begin{align*}
\text{MSE} = & \ (\text{MSE}_{u_{in}} + \text{MSE}_{((-\Delta)^{-1}u)_{in}}) + (\text{MSE}_{u_{b}} + \text{MSE}_{((-\Delta)^{-1}u)_{b}}) + \text{MSE}_{F} \\& +\text{MSE}_{Lap} + \text{MSE}_{v} + \text{MSE}_{zm}.
\end{align*}
Moreover, since certain components have more importance over the others, such as the mean square error on the initial condition for its ground truth availability, adjustable weights are applied to each component:
\begin{align*}
\text{MSE} = & \  {W_{u_{in}}}\text{MSE}_{u_{in}} +{W_{((-\Delta)^{-1}u)_{in}}}\text{MSE}_{((-\Delta)^{-1}u)_{in}}  \\ & + {W_{u_{b}}}\text{MSE}_{u_{b}} + {W_{((-\Delta)^{-1}u)_{b}}}\text{MSE}_{((-\Delta)^{-1}u)_{b}} \\ & + {W_{F}}\text{MSE}_{F} +{W_{Lap}}\text{MSE}_{Lap} +{W_{int}}\text{MSE}_{int}+{W_{zm}}\text{MSE}_{zm}
\end{align*}


\newpage
\section{Numerical Results}\label{Results}
Now we have all the pieces regarding how to apply PINNs on ACOK. We can combine them and implement them in Python. The code is based on the work of Maziar Raissi et al. \cite{PINNs}.
We now show the results of the modified PINNs on 1-dimensional ACOK with different time intervals, ranging from  $1\times 10^{-3}$s to $1\times 10^{-2}$s. Fine-tuning on hyper-parameter is needed to reach the best results, but it is out of our current research scope. 

\begin{figure}[H]
\centering
\subfigure[]{\includegraphics[scale=0.475]{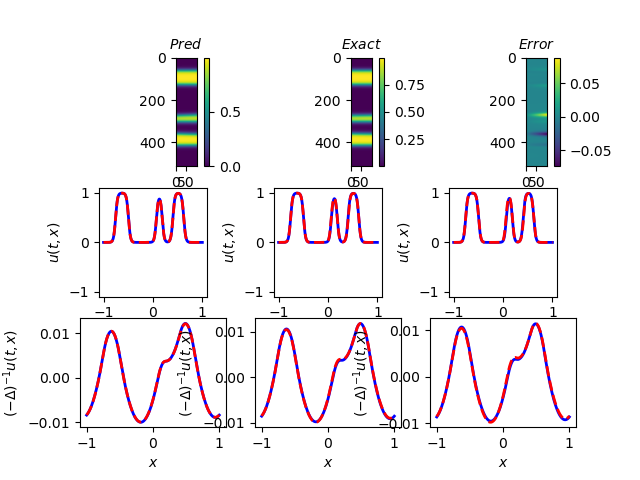}}
\subfigure[]{\includegraphics[scale=0.475]{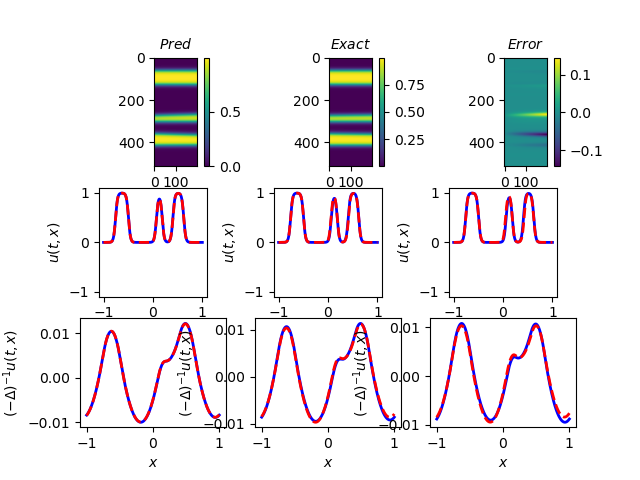}}

\subfigure[]{\includegraphics[scale=0.475]{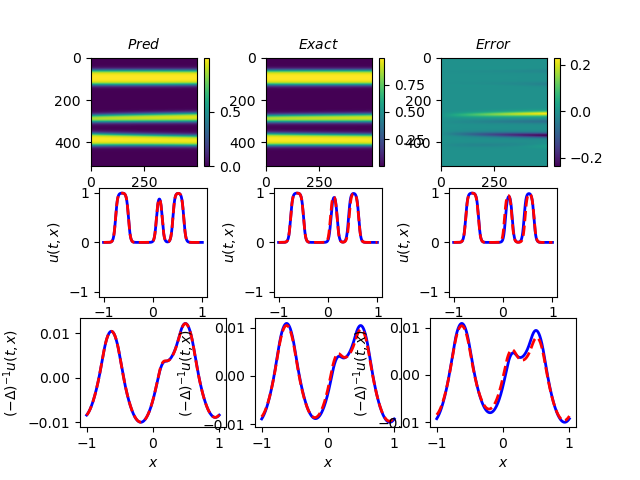}}
\subfigure[]{\includegraphics[scale=0.475]{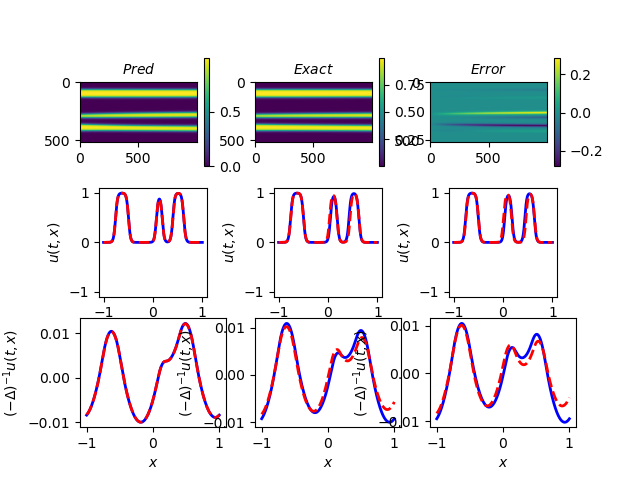}}
\caption{Numerical approximations of the ACOK equation using the modified PINNs with various setups. (a)$t = 1\times 10^{-3}$, $N_0 = 500$, $N_b = 95$, $N_f = 2\times 10^{4}$, $N_{t_{uni}} = 20$; (b) $t = 2\times 10^{-3}$, $N_0 = 500$, $N_b = 195$, $N_f = 4\times 10^{4}$, $N_{t_{uni}} = 40$; (c) $t = 5\times 10^{-3}$, $N_0 = 500$, $N_b = 495$, $N_f = 1\times 10^{5}$, $N_{t_{uni}} = 100$; (d) $t = 1\times 10^{-2}$, $N_0 = 500$, $N_b = 995$, $N_f = 2\times 10^{5}$, $N_{t_{uni}} = 200$. }
\label{fig:example1}
\end{figure}

\Cref{fig:example1} are run for $495$ epochs, with weights $W_{u_{in}} = 1\times 10^5$, $W_{((-\Delta)^{-1}u)_{in}} = 5 \times 10^6$, $W_{u_{b}} = 1$, $W_{((-\Delta)^{-1}u)_{b}} = 30$, $W_{F} = 1$, $W_{Lap} = 500$, $W_{int} = 1$, $W_{zm} = 30$. There are $10$ hidden layers in $\text{Net}_u$, each containing 20 neurons, while only $3$ hidden layers with $10$ neurons for each layer were needed for $\text{Net}_v$, since it is a much simpler neural network.

We can see that the modified PINNs perform remarkably well in smaller time periods. However, as the time interval increases, the prediction of $(-\Delta)^{-1}u$ at the right tail is getting notably worse. 
There are several ways to address this problem. One approach is the Time-Adaptive Method \cite{timeadap}, or similarly, the sequence-to-sequence learning \cite{FailurePINNs}, which learns only one small-time period at a time and then uses the prediction as to the initial data for the next period. Since our model approximates well with short time intervals, we can compile all the sub-results together to get the entire solution. For example, we could use the data of the last time column in \Cref{fig:example1}(a) as the initial data for the next 100 columns and repeat the process in 100 column increments to achieve good results.

\begin{figure}[H]
\centering
\subfigure[]{\includegraphics[scale=0.175]{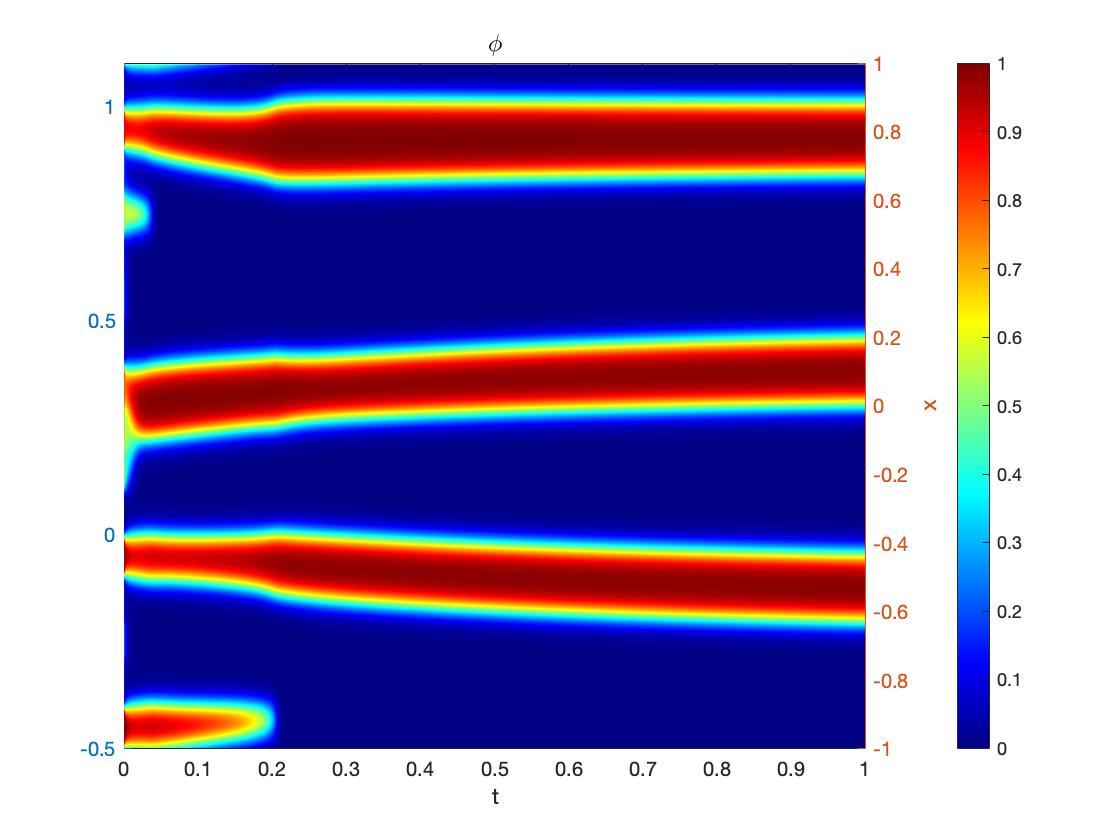}}
\subfigure[]{\includegraphics[scale=0.35]{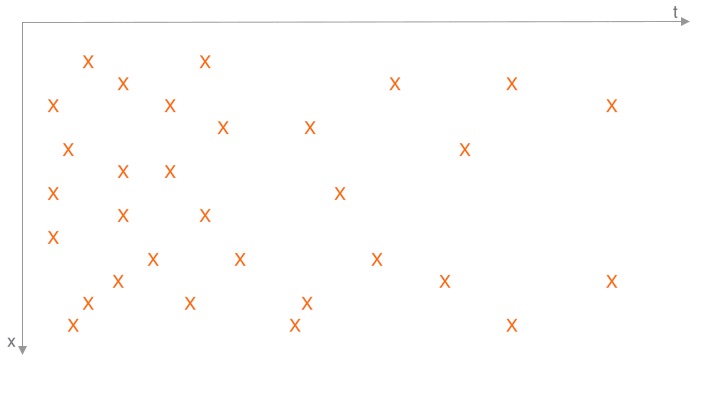}}
\caption{An illustration on how to improve the modified PINNs. (a) numerical results of the 1D ACOK equation with spectrum method; (b) scaling sample points. }\label{Fig:Matlab}
\end{figure}

Furthermore, the problem can be improved by rescaling the sampling points. As shown in \Cref{Fig:Matlab}(a), changes in phase separation occur most rapidly during the earlier time period. Therefore, we propose rescaling the sampling points to be denser in earlier time columns. 
Recall \Cref{Fig:Uniform2}, we continue to use the LHS sampling, but rescale the points polynomially to achieve the denser in the beginning effect, as in  \Cref{Fig:Matlab}(b).

\begin{figure}[H]
\centering
\subfigure[Original LHS sampling.]{\includegraphics[scale=0.475]{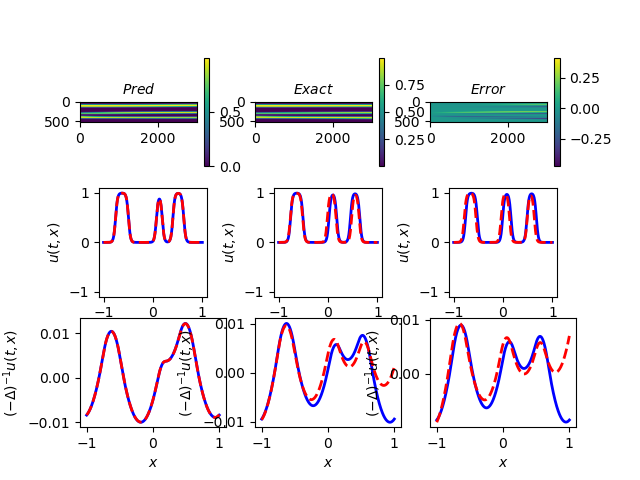}}
\subfigure[$x^2$ rescaling sampling.]{\includegraphics[scale=0.475]{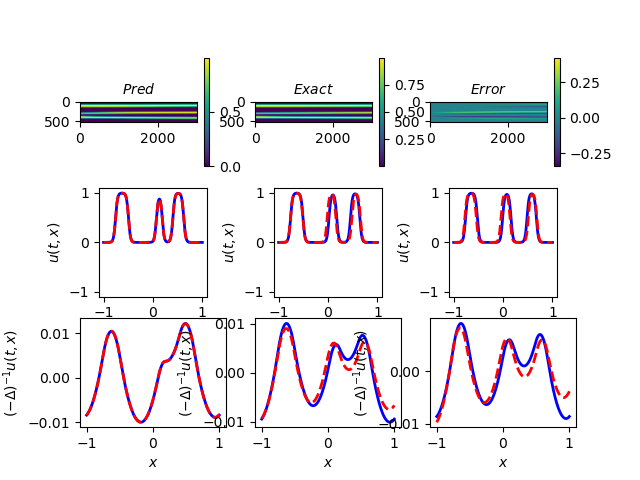}}
\caption{A comparison of different sampling strategy. (a) shows the original LHS sampling; (b) shows the $x^2$ rescaling sampling. In both cases, we use $t = 3\times 10^{-2}$, $N_0 = 500$, $N_b = 2995$, $N_f = 6\times 10^{5}$, and $N_{t_{uni}} = 10$. }\label{Fig:Regular}
\end{figure}

The results are shown as follows, \Cref{Fig:Regular}(a) is the result with original LHS sampling, with time interval being $3\times 10^{-3}$s.
Notice that the right tail on $(-\Delta)^{-1}u$ fits quite badly in this case due to the long time interval. \Cref{Fig:Regular}(b) is the result with a $x^2$ rescaling on the points. Although it does not fix the problem entirely, we can see that the tail on $(-\Delta)^{-1}u$ fits much better with a $x^2$ rescaling. In some cases, $x^2$ rescaling does not improve the result much. Therefore, $x^3$ rescaling could potentially be used to improve the fitting further.

\section{Discussion and future work} \label{conclusion}

During implementation, we have discovered that our modified PINNs are challenging to train due to their many hyper-parameters. Some hyper-parameters, such as learning rate, number of hidden layers, and number of neurons of each layer, are fixed during our training. However, the weights, number of epochs, and number of sampling points need fine-tuning to achieve the best results.

The most challenging part of our fine-tuning is finding a well-balanced set of weights. Each weight controls a different part of the estimation. For instance, $W_{u_{in}}$, $W_{((-\Delta)^{-1}u)_{in}}$ controls the fitting of initial data, which is undoubtedly an essential part of our loss function since the initial data is the only ground truth data we have. But they differ from each other, as we later find out that the scale of $((-\Delta)^{-1}u)_{in}$ is 100 times smaller than that of $u_{in}$. Therefore, we also consider that by increasing $W_{((-\Delta)^{-1}u)_{in}}$ to counter the scale issue. A similar scaling technique is applied to $W_{u_{b}}$ and $W_{((-\Delta)^{-1}u)_{b}}$, which controls the periodic boundary conditions. 

$W_{F}$ dominates the physics information in our model, while $W_{Lap}$ controls the accuracy of the approximation of the long-range interaction term. We have learned that the accuracy of the long-range interaction term is one of the key players in our approximation. $W_{int}$ directs the volume constraint and, therefore, can be increased if the predicted result has an incorrect volume. $W_{zm}$ applies only to the $((-\Delta)^{-1}u)_{in}$, which means the scaling problem also needs to be taken into consideration. We have also discovered that the interval of a number of epochs is quite important. The model will not train well with too few epochs. However, if we boost the number of epochs to the scale of the thousands, the model is over-trained and does not perform well.

Moreover, we tune the number of the sampling points using the principle that we want to keep the randomness while at the same time preserving as much of the known information as possible. We also remain attentive to the efficiency of the program. Again, the initial data, being the only piece of ground truth information available, need to have the most sampling points proportionally. We do not take all the initial data because we want to keep the randomness. We also sample a relatively large percentage of the boundary points for the same reason. As the total number of the grid points is quite large for the internal points, we cannot take as many as the sampling points while keeping the program's efficiency. Randomness also plays a crucial part here. Therefore, we only select less than half of the points. The number of randomly selected $t$ columns in the uniform mesh, $N_{t_{uni}}$, follows the same principle.

To increase the randomness, we also implemented a mini-batch method. The idea is that instead of considering all the points simultaneously, we consider a small subset of the points at a time. In theory, it should need fewer epochs and converge faster. However, as the computation in our case is quite complex, we have discovered that the mini-batch method takes significantly longer time than the regular approach. 

We are currently working on the 2-dimensional case of this problem. As the number of points has drastically increased for the same time period, we have run into technical issues, in particular, insufficient memory capacity. To address this issue, we have to dramatically decrease the sampling points, even on the initial data. Also, we are keeping the time interval relatively small for the same reason. Additionally, for the boundary condition, we now have the upper and lower bound and the front and back bound. A more complicated back-propagation further lengthens the training time due to the second-order derivative. These are the difficulties we have encountered in the 2-dimensional case.

Possible future work includes investigating the inverse problem for nonlocal phase-field models with unknown nonlocal operators. We could train our model to learn the dynamics using the existing data and obtain a set of trained weights and bias representing the unknown operator to compare with a number of known possible operators. Furthermore, we could also explore other methods to improve the efficiency of the training process. The current one-dimensional results run within an hour. However, as the structure of the DNN becomes more complicated, training time increases drastically, likely caused by the more complex computations and the increased number of epochs. One possible solution is to start the currently randomized initial weights and bias closer to the solution, which will significantly lower the number of epochs, and by extension, decrease the training time overall. Additionally, various other machine learning methods, such as FNO, conventional neural networks, and transformer-based methods, can be applied to PDEs. These are also some exciting directions to explore.

\section*{Acknowledgements}
The work of J. Xu and Y. Zhao is supported by a grant from the Simons Foundation through Grant No. 357963 and NSF grant DMS-2142500. J. Zhao would like to acknowledge the support from National Science Foundation, United States, with grant DMS-2111479.

  
  


\addcontentsline{toc}{section}{Bibliography}

\bibliographystyle{plain} 
\bibliography{references} 


\end{document}